\documentclass{article}

\usepackage{amsmath}
\usepackage{amssymb}
\usepackage{color}

\usepackage{latexsym}
\usepackage[english]{babel}
\usepackage[dvips]{graphicx}
\usepackage{verbatim}
\newcommand{\Gq}{{\gg}_q}

\renewcommand{\ll}{\langle}
\newcommand{\rr}{\rangle}

\newcommand{\Zz}{\mathcal{Z}}

\renewcommand{\gg}{{\bf G}}

\newcommand{\bbibitem}{\bibitem}

\newcommand{\llabel}[1]{{\label{#1}}}

\renewcommand{\r}[1]{(\ref{#1})}

\newcommand{\ex}[1]{} 

\font\tenmsb=msbm10
\font\sevenmsb=msbm7
\font\fivemsb=msbm5
\newfam\msbfam
\textfont\msbfam=\tenmsb
\scriptfont\msbfam=\sevenmsb
\scriptscriptfont\msbfam=\fivemsb

\newcommand{\bi}{\begin{itemize}}
\newcommand{\ei}{\end{itemize}}
\newcommand{\bd}{\begin{description}}
\newcommand{\ed}{\end{description}}

\newcommand{\bqn}{\begin{eqnarray}}
\newcommand{\eqn}{\end{eqnarray}}
\newcommand{\eqnn}{\nonumber\end{eqnarray}}
\newcommand{\eqnl}[1]{\llabel{#1}\end{eqnarray}}

\newcommand{\nn}{\nonumber}
\newcommand{\ba}[1]{\begin{array}{#1}}
\newcommand{\ea}{\end{array}}

\newcommand{\R}
{\mathbb{R}}

\newcommand{\N}{\mathbb{N}}

\newcommand{\fine}{\end{document}}

\def \trait (#1) (#2) (#3){\vrule width #1pt height #2pt depth #3pt}
\def \qed{\hfill
        \trait (0.1) (6) (0)
        \trait (6) (0.1) (0)
        \kern-6pt   
        \trait (6) (6) (-5.9)
        \trait (0.1) (6) (0)
\medskip}
\def \qedmio{\hfill
             \trait (8) (8) (-0.1)
             \medskip}
\def \quadp{{\Huge $\qedmio$}}

\newtheorem{Theorem}{\bf Theorem}
\newtheorem{ml}[Theorem]{\bf Lemma}
\newtheorem{mo}{\bf \underline{{\sl Observation}}}
\newtheorem{mcc}[Theorem]{\bf Corollary}
\newtheorem{Definition}[Theorem]{\bf Definition}
\newtheorem{mpr}[Theorem]{\bf Proposition}
\newtheorem{mproperty}[Theorem]{\bf Property}

\newtheorem{mrem}[Theorem]{\bf \underline{{\sl Remark}}}

\newcommand{\bt}{\begin{Theorem}}
\newcommand{\et}{\end{Theorem}}
\newcommand{\bl}{\begin{ml}}
\newcommand{\el}{\end{ml}}
\newcommand{\bo}{\noindent\begin{mo}\rm}
\newcommand{\eo}{\end{mo}}
\newcommand{\bp}{\begin{mpr}}
\newcommand{\ep}{\end{mpr}}
\newcommand{\bc}{\begin{mcc}}
\newcommand{\ec}{\end{mcc}}
\newcommand{\bdeff}{\begin{Definition}}
\newcommand{\edeff}{\end{Definition}}
\newcommand{\bproperty}{\begin{mproperty}}
\newcommand{\eproperty}{\end{mproperty}}
\newcommand{\brem}{\begin{mrem}\rm}
\newcommand{\erem}{\end{mrem}}

\newcommand{\bfi}{\begin{figure}}
\newcommand{\efi}{\end{figure}}

\newcommand{\ppotR}[3]
{

\begin{figure}\begin{center}
~\includegraphics[width=#3truecm]{#1.eps}\\
\caption{#2}
\llabel{#1}
\end{center}
\end{figure}
\noindent$\!\!$}

\newcommand{\lam}{\lambda}

\newcommand{\g}{\gamma}
\newcommand{\al}{\alpha}
\newcommand{\eps}{\varepsilon}

\newcommand{\om}{\omega}
\newcommand{\Om}{\Omega}

\newcommand{\con}{{\cal C}}

\newcommand{\rvd}{rank-varying distribution}

\newcommand{\ar}{rank-varying sub-Riemannian structure}
\newcommand{\Ar}{Rank-varying sub-Riemannian structure}
\newcommand{\VecM}{\mathrm{Vec}(M)}
\newcommand{\bD}{\Delta}

\newcommand{\proof}{{\bf Proof. }}

\newcommand{\brs}{\begin{eqnarray*}}
\newcommand{\ers}{\end{eqnarray*}}
\newcommand{\br}{\begin{eqnarray}}
\newcommand{\er}{\end{eqnarray}}

\newcommand{\rw}{\rightarrow}

\newcommand{\lp}{\left(}
\newcommand{\rp}{\right)}
\newcommand{\mt}{\mapsto}

\newcommand{\be}{\begin{equation}}
\newcommand{\ee}{\end{equation}}

\begin{document}
\begin{center} \noindent
{\LARGE{\sl{\bf A Gauss-Bonnet-like Formula on Two-Dimensional 

\vskip .2cm
Almost-Riemannian 
Manifolds}}}\footnote{The second and third authors have been supported by the ``Research in Pairs'' program of the MFO 
(Mathematisches Forschungsinstitut Oberwolfach).}
\end{center}

\vskip 1cm
\begin{center}
Andrei A. Agrachev, Ugo Boscain

{\footnotesize SISSA-ISAS,
Via Beirut 2-4, 34014 Trieste, Italy
 \\{\tt agrachev(-at-)sissa.it, boscain(-at-)sissa.it}}

\vspace{4mm}

Mario Sigalotti

{\footnotesize
Institut \'Elie Cartan, UMR 7502 Nancy-Universit\'e/CNRS/INRIA,
POB 239, 54506 Vand\oe uvre-l\`es-Nancy,
France {\tt Mario.Sigalotti(-at-)inria.fr}}

\end{center}

\begin{center}
\today
\end{center}
\vspace{1cm}

\begin{quotation}\noindent  {\bf\em Abstract --- 
We consider a generalization of Riemannian geometry that naturally arises
in the framework of control theory. Let $X$ and $Y$ be 
two smooth vector fields on a two-dimensional  manifold $M$. If $X$ and 
$Y$ are everywhere linearly independent, then  they define a classical 
Riemannian metric on $M$ (the metric for which they are orthonormal) and 
they give to $M$ the structure of metric space.
If $X$ and $Y$ become linearly dependent somewhere on $M$, 
then the corresponding Riemannian metric has singularities, but under 
generic conditions the metric structure is still well defined.  
Metric structures that can be defined locally in this way 
are called almost-Riemannian structures.  They are special cases of \ar s,
which are naturally defined in terms of submodules of the space of smooth vector fields on $M$. 
Almost-Riemannian structures show
interesting phenomena, in particular for what concerns the 
relation between curvature, presence of conjugate points, and topology of 
the manifold.
The main result of the paper is a generalization to almost-Riemannian 
structures of the Gauss-Bonnet formula.

}\end{quotation}


\vskip 1cm
\noindent
{\bf MSC-class:} 49j15, 53c17

\vskip 2cm
\begin{center}
PREPRINT SISSA  55/2006/M
\end{center}

\newpage


\section{Introduction}

Let $M$ be a two-dimensional smooth manifold and consider a pair of smooth vector fields $X$ and $Y$ on $M$. 
If the pair $X$, $Y$ is Lie bracket generating, i.e., if
$\mathrm{span}\{X(q),Y(q),[X,Y](q),[X,[X,Y]](q),\ldots\}$ is 
full-dimensional at every $q\in M$, then the
control system
\bqn
\label{ff}
\dot q=u X(q)+v Y(q)\,,~~~u^2+v^2\leq 1\,,~~~q\in M\,,
\eqn
is completely controllable and  the 
minimum-time function defines a 
continuous distance $d$ on $M$. 
When $X$ and $Y$ are everywhere linear independent 
(the only possibility for this to happen is that $M$ is parallelizable), 
such distance is Riemannian and it 
corresponds 
to the metric for which $(X,Y)$ is an orthonormal moving frame. 
The idea is to study the geometry obtained starting from a pair of vector fields which 
may  become collinear. 
Under generic hypotheses, the set $\Zz$ (called {\it singular locus}) of 
points of $M$ at which  $X$ and $Y$ are parallel 
is  a one-dimensional embedded submanifold of $M$ (possibly disconnected). 

Metric structures
that can be defined {\it locally} by a pair of vector fields 
$(X,Y)$ through \r{ff} 
are called almost-Riemannian structures.

An almost-Riemannian structure can equivalently be seen as a locally finitely generated Lie bracket generating
$\con^\infty(M)$-submodule $\bD$ of $\VecM$, the space of smooth vector fields on $M$, 
endowed with a bilinear, symmetric map $G:\bD\times\bD\to \con^\infty(M)$ which is positive definite (in a suitable sense). 
A pair of vector fields $X$, $Y$ in $\bD$ is said to be orthonormal on some open set $\Omega$ if $G(X,Y)(q)=0$ and 
$G(X,X)(q)=G(Y,Y)(q)=1$ for every $q\in\Omega$.

An almost-Riemannian structure is said to be {\it orientable} if there exists a {\it volume 
form}, i.e., a bilinear, skew-symmetric, 
non-degenerate form
$\omega:\bD\times\bD\to \con^\infty(M)$. In this case it is possible to normalize 
$\omega$ in such that $|\omega(X,Y)|=1$ on $\Omega$ 
for every open subset $\Omega$ of $M$ and every local orthonormal frame $(X,Y)$ on $\Omega$.

It is interesting to notice that it is possible to define non-orientable 
almost-Riemannian structures on orientable manifolds and orientable 
almost-Riemannian structures on non-orientable manifolds.

We say that an almost-Riemannian structure is {\it trivializable} if $\Delta$ is 
 globally generated by a pair of vector fields defined on $M$. Trivializable 
almost-Riemannian structures are always orientable.

The singular locus $\Zz$
can be defined on $M$ 
as the set where the linear subspace $\Delta(q)=\{V(q)\mid V\in\bD\}$ 
of $T_q M$
is not full-rank. 
An almost-Riemannian structure is Riemannian if and only if $\Zz=\emptyset$.

A famous 
example of 
genuine
almost-Riemannian structure 
is provided by 
the 
Grushin plane, 
which is the almost-Riemannian structure on $\R^2$ 
for which the vector fields $X(x,y)=(1,0)$ and
$Y(x,y)=(0,x)$ form a pair of orthonormal generators. 
(See Section \ref{ss-grushin} and \cite{bellaiche,algeria}. 
The model was originally introduced in the context of
hypoelliptic operator theory 
\cite{FL1,FL2,grusin1,grusin2}.)
Notice that the 
singular locus
is indeed nonempty, being equal to
the $y$-axis. 
Another example of (trivializable) almost-Riemannian structure
has appeared in problems of control of quantum mechanical systems 
(see \cite{q1,q4}). In this case $M=S^2$ represents a suitable state space reduction of a 
three-level quantum 
system while the orthonormal generators $X$ and $Y$ are two infinitesimal rotations along two orthogonal axes,  
modeling 
the action on the system of
two lasers in the rotating wave approximation (see Section \ref{s-sketch}). 

Such examples, and the naturalness 
of the construction leading to the 
definition of almost-Riemannian structure, motivate the study of 
general properties of such geometry, which exhibits 
many interesting
features. 
One can check, for instance, that even in the case where the Gaussian curvature is everywhere negative 
(where it is defined, i.e., on $M\setminus\Zz$) 
geodesics may have conjugate points. For this reason it 
seems interesting to 
analyze the relations between the curvature, the presence of conjugate points, 
and the topology of the manifold (see also  \cite{agra-gauss}). 
After providing a characterization of generic almost-Riemannian structures
by means of local normal forms, in this paper
we start this program by proving 
a generalization of the Gauss-Bonnet formula. 
Let $M$ be compact and oriented, and endow it with an orientable almost-Riemannian structure. 
Denote by 
$K:M\setminus\Zz\rw\R$ 
the Gaussian curvature. 
The first difficulty in order to extend the Gauss-Bonnet formula is to give a meaning 
to $\int_M K\, d A$, 
the integral of $K$ on $M$ with respect to the Riemannian density $d A$ induced by the 
Riemannian metric on $M\setminus\Zz$. 
In the two examples cited above one can check, for instance, that, as $q$ 
approaches $\Zz$,
$d A$ diverges, while 
$K(q)$, which is everywhere negative, tends to  
$-\infty$. 

The idea is to replace $K\,dA$ with a signed version of it. 
A natural choice is $K\,dA_s$, where $dA_s$ is a volume form intrinsically associated 
with the almost-Riemannian structure  on $M\setminus\Zz$.

Our goal is to prove the existence and to assign a value to the limit
\bqn
\label{limit}
\lim_{\eps\searrow0}\int_{\{q\in M \mid	d(q,\Zz)>\eps\}}K(q) dA_s,
\eqn
where
$d(\cdot,\cdot)$ is the distance globally defined by the almost-Riemannian structure on $M$.

The goal will be attained under the following additional assumption.
Generically the singular locus $\Zz$ is smooth and $\Delta(q)$ is one-dimensional at every 
point of $\Zz$. 
We say that 
$q\in\Zz$ is a {\it tangency point} if 
$\Delta(q)$ is tangent to $\Zz$.  
Under generic assumptions, $\Zz$ contains only a discrete set of tangency 
points.
 The hypothesis under which the main results of the paper are obtained is that
 $\Zz$ contains no such point. 

Define $M^+$ (respectively, $M^-$) as the subset of $M\setminus\Zz$ 
on which the orientation defined by $dA_s$ coincides with (respectively, is opposite to) that of $M$. 
If $M$ has no tangency point, 
then the limit \r{limit} turns out to exist and is equal to 
 $2\pi (\chi(M^+)-\chi(M^-))$, where $\chi$ denotes the Euler characteristic.

When the almost-Riemannian structure 
 is trivializable,
 we have that $\chi(M^+)=\chi(M^-)$ and thus the 
 limit \r{limit}
is equal to zero. Once applied to the special subclass of Riemannian structures, such result simply states that 
the  integral of the curvature of a parallelizable compact oriented surface (i.e., the torus) is equal to zero. 
In a sense,
in the standard Riemannian construction 
the topology of the surface gives a constraint  
on the total curvature
through the Gauss-Bonnet formula,
whereas for 
an almost-Riemannian structure induced by a single pair of vector fields 
the total curvature is equal to zero and the topology of the manifold 
constrains the metric to be singular on a suitable set.

It is interesting to notice that every oriented compact surface 
can be endowed with a trivializable almost-Riemannian structure 
satisfying the requirement that there are no
tangency points.

The paper is organized as follows. 
In Section~\ref{s-definitions}, we introduce two equivalent definitions 
of \ar\ on a manifold of any dimension, first by using the language of moduli 
and then by identifying it with an atlas of orthonormal frames. 
Rank-varying sub-Riemannian structures have already been studied, from a 
different perspective, in \cite{jean1,jean2}.

A notion of orientability for \ar\ is then introduced.
Almost-Riemannian structures are defined as 
\ar s of maximal rank.

Starting from Section~\ref{s-dim2} we focus on the case of almost-Riemannian structures 
on two-dimensional manifolds.
Geodesics associated with such structures are characterized
in Section \ref{s-minimizers} using  
the 
Pontryagin Maximum Principle. 
In Section \ref{ss-grushin} we study the Grushin 
plane, for which we compute the cut and the conjugate loci.

In Section \ref{s-generic} we provide local normal forms
for generic almost-Riemannian structures, which are used in Section \ref{s-main} to prove a generalization
of the  Gauss-Bonnet formula to almost-Riemannian structures without tangency points.
The formula is then specialized to the case of trivializable 
almost-Riemannian structures. 
In Section \ref{s-sketch} we show that every compact orientable two dimensional manifold admits 
a trivializable  almost-Riemannian structure with no tangency points.



\section{Rank-varying distributions and sub-Riemannian structures}\label{s-definitions}

Let $M$ be a   $n$-dimensional smooth manifold.
Recall that $\mathrm{Vec}(M)$, the set of 
smooth vector fields on  $M$, 
is naturally endowed with the structure of $\con^\infty(M)$-module.
Given an open subset $\Omega$ of $M$, a submodule 
$\Delta$ of $\VecM$ is said to be {\it generated on $\Omega$}  
by the vector fields $\{V_1,\ldots,V_m\}$ if every $V\in \bD$ 
can be written as
$V=a_1 V_1+\cdots +a_m V_m$ on $\Omega$ where $a_1,\dots,a_m$ belong to $\con^\infty(M)$.

\bdeff
A {\it $(n,k)$-\rvd} is 
a pair $(M,\bD)$ where 
$M$ is a $n$-dimensional smooth manifold, $\bD$ is a  
submodule  of $\mathrm{Vec}(M)$, and $k\leq n$ 
is such that 
for every $q\in M$ and every small enough 
neighborhood $\Omega^q$ of $q$, the restriction to $\Omega^q$ of $\bD$ is 
generated by $k$ vector fields and cannot be generated by less than $k$ 
vector fields. 
\edeff


From now on the expression $\bD(q)$ will denote the linear subspace $\{V(q)\mid  V\in \bD\}\subset  T_q M$.
When the dimension of $\bD(q)$ is independent of $q$, we recover the standard 
definition of distribution as a smooth field of linear subspaces of $T_qM$.
Notice that $\bD$ cannot be identified with the map $q\mapsto\Delta(q)$. Indeed, it can happen that two different 
moduli $\bD_1$ and $\bD_2$ are such that 
$\bD_1(q)=\bD_2(q)$ for every $q\in M$. Take for instance $M=\R$ and  $\bD_1$,  $\bD_2$ generated,  respectively, by 
$F_1(x)=x$, $F_2(x)=x^2$.

Denote by  $\mathrm{Lie}(\bD)$  the smallest Lie subalgebra of
$\mathrm{Vec}(M)$
containing $\bD$ and let $\mathrm{Lie}_q(\bD)=\{V(q)\mid V\in \mathrm{Lie}(\bD)\}$ for
every $q\in M$.
We say that $(M,\bD)$ satisfies the Lie bracket generating condition if
$\mathrm{Lie}_q(\bD)=T_q M$ for every $q\in M$.
We also introduce the flag of a \rvd\ $(M,\bD)$ as the sequence of submodules
$\bD_0=\bD\subset \bD_1\subset\ldots \subset\bD_m \subset \cdots$ 
defined through the recursive formula
\bqn
\bD_{k+1}=\bD_k+[\bD,\bD_k].
\eqnl{flag}
As above, we let $\bD_{m}(q)=\{V(q)\mid V\in \bD_m\}$.

In order to provide an example of \rvd, let us introduce 
the {\it Grushin distribution}.
Take as $M$ the plane $\R^2$ and let $\bD$ be generated by the vector 
fields $F_1(x,y)=(1,0)$ 
and $F_2(x,y)=(0,x)$. Then $(\R^2,\bD)$ is a $(2,2)$-\rvd.
Notice that $\Delta(q)$ is equal to $\R\times\{0\}$ when $q$ is on the $y$-axis and to $\R^2$ elsewhere. The Grushin distribution is Lie bracket generating since $\Delta_2(q)=\R^2$ for every $q\in \R^2$. \\\\

Crucial in what follows is the notion of {\it 
generic} $(n,k)$-\rvd. 
Denote by ${\cal W}$ the  $\con^2$-Whitney topology defined on $\mathrm{Vec}(M)$ 
and by $(\mathrm{Vec}(M),{\cal W})^k$ the product of $k$ copies of  
$\mathrm{Vec}(M)$ endowed with the corresponding  product  topology. 
We recall that if $M$ is compact (as it is the case in most of what follows), then ${\cal W}$ 
is the standard $\con^2$ topology.
\bdeff
A property $(P)$ defined for $(n,k)$-\rvd s is said to be {\it generic}
if there exists an  open and dense  subset ${\cal O}$ of $(\mathrm{Vec}(M),{\cal W})^k$ such 
that $(P)$ holds for every $(n,k)$-\rvd\ which is 
generated by elements of  ${\cal O}$.
\edeff
E.g., generically, a $(n,k)$-\rvd\ is  Lie bracket generating.
\subsection{Orientable \rvd s}
Let $(M,\bD)$ be a  $(n,k)$-\rvd. 
A  $k$-form on  $(M,\bD)$
is a multilinear 
skew-symmetric map
$$\omega:\underbrace{\bD\times\cdots\times\bD}_{k~\mathrm{ times}}\to \con^\infty(M).$$
We say that a $k$-form $\omega$ is 
a {\it volume form} 
if, 
for every $q\in M$, there exist $k$ vector fields $F_1,\dots,F_k\in\bD$ such that 
$\omega(F_1,\dots,F_k)(q)\neq0$.
\bdeff
We say that a $(n,k)$-\rvd\ $(M,\bD)$ is orientable if it admits 
a volume form,
otherwise we say that $(M,\bD)$ 
is non-orientable.  
\edeff
Notice that a \rvd\ can be orientable even if $M$ is a non-orientable manifold (see example below).  
However, the distribution $(M,\mathrm{Vec}(M))$ is orientable if and only if $M$ is. 
\brem
Thanks to its multilinearity, 
a volume form 
is completely 
characterized by its action on the generators.
Given a trivializable  \rvd\ $\bD$ and a  global system 
of generators  $F_1,\ldots, F_k$, 
the equality $\omega(F_1,\ldots,F_k)=1$
uniquely defines 
a volume form 
on $\bD$. Hence every  trivializable 
\rvd\ is orientable. 
\erem
\brem\label{r-tensor}
Let $\om$ be a $k$-form on a $(n,k)$-\rvd\ $(M,\bD)$. 
Then $\om$ acts as a tensor on the open subset of $M$ made of points $q$ such that
the dimension of $\bD(q)$ is equal to $k$, i.e., for every $V_1,\dots,V_k\in \bD$ the value of $\om(V_1,\dots,V_k)(q)$ depends only 
on $V_1(q),\dots,V_k(q)$. 
Indeed, let $\{F_1,\dots,F_k\}$ be a local system of generators of $\bD$ on a neighborhood $\Om$ of $q$ and take $a_{ij}
\in\con^\infty(M)$, $1\leq i,j\leq k$, such that $V_i=\sum_{j=1}^k a_{ij}F_j$ on $\Om$.
The multilinearity and skew-symmetricity of $\om$ implie that 
$\om(V_1,\dots,V_k)(q)=\mathrm{det}(a_{ij}(q)) \om(F_1,\dots,F_k)(q)$. Therefore, 
$\om(V_1,\dots,V_k)(q)$ depends on $V_1,\dots,V_k$ only through 
the matrix $(a_{ij}(q))$, which is 
uniquely determined by $V_1(q),\ldots,V_k(q)$.
\erem

Let us present some example of orientable and 
non-orientable \rvd s. All these examples are  $(2,2)$-\rvd s, since this  
is our main case of interest in the following.\\\\
{\bf The Grushin distribution.}
Let $M=\R^2$ and recall that $\bD$ is generated by the vector 
fields $F_1(x,y)=(1,0)$ 
and $F_2(x,y)=(0,x)$.
The distribution $\Delta$ is orientable. A  
volume form 
can be  
defined by 
\bqn
\omega(V_1,V_2)(x,y)=
\lim_{x'\to x}\frac1{x'}dx\wedge dy~(V_1(x',y),V_2(x',y)),
\eqnn
for every pair of vector fields $V_1$, $V_2$ belonging to $\bD$. 
Equivalently we could have defined $\omega$ on the generators $F_1$, 
$F_2$ 
by requiring that $\omega(F_1,F_2)(x,y)=1$.

Notice that $\omega$ is not a  tensor on the $y$-axis.\\\\
{\bf A non-orientable \rvd\ on the torus.}
Let $M$ be the two-dimensional torus $[-\pi,\pi]\times[-\pi,\pi]$ with 
the 
standard identifications. Consider 
the open covering of $M$ given by
\bqn
\Omega^1=(-\pi/2,\pi/2)\times[-\pi,\pi],~~
\Omega^2=
\big([-\pi,-\pi/4)
\cup(\pi/4,\pi]
\big)
\times
[-\pi,\pi].\nn
\eqn
Let $\bD$ be generated by the vector fields 
\bqn
\ba{lll}
F_1^1=(1,0),&
F_2^1=(0,\sin x),&\mbox{ on }\Omega^1\\
F_1^2=(1,0),&
F_2^2=(0,1),&\mbox{ on }\Omega^2.\ea
\eqnn
This \rvd\ is non-orientable. 
In fact a volume form 
$\omega$ should acts on  the local generators 
as
\bqn
\omega(F^1_1,F^1_2)(q)=f_1(q)\mbox{ for every 
}q\in\Omega^1,~~~~\omega(F^2_1,F^2_2)(q)=f_2(q)\mbox{ for every
}q\in\Omega^2,\nn\eqn
where $f_1$ and $f_2$ are two  never-vanishing smooth functions. 
On $\Omega^1\cap\Omega^2$ we would have 
\bqn
f_1(x,y)=\omega(F_1^1,F_2^1)(x,y)=\sin (x)\omega(F_1^2,F_2^2)(x,y)=\sin (x) f_2(x,y),
\eqnn
contradicting the constant-sign assumption on $f_1$, $f_2$.
As a consequence $\Delta$ is not trivializable.\\\\
{\bf An orientable \rvd\ on the Klein bottle.} 
Let $M$ be the Klein bottle seen as the square 
$[-\pi,\pi]\times[-\pi,\pi]$ with 
the identifications $(x,-\pi)\sim(x,\pi)$,  
$(-\pi,y)\sim(\pi,-y)$. 
Consider the vector fields 
\bqn
F_1(x,y)=(1,0),~~F_2(x,y)=(0,\sin(2x)),
\eqnn
which are well defined on $M$.
The distribution generated by $F_1$ and $F_2$ is orientable since it is 
trivializable.\\\\

\subsection{\Ar s}
In this section we see how to introduce a a smoothly-varying Riemannian structure on every subspace $\bD(q)$.
\bdeff  
\label{d-ars-1}
A $(n,k)$-\ar\  is a triple ${\cal S}=(M,\bD,G)$, where $(M,\bD)$ is a Lie bracket generating 
$(n,k)$-\rvd\ and 
$G:\bD\times\bD\to\con^\infty(M)$ 
is a  symmetric,  positive definite 
bilinear map, i.e., a map such that for every $V,W\in\bD$ and  
$f\in\con^\infty(M)$ we have
\bqn
&&G(V,W)=G(W,V),\nn\\
&&G(fV,W)=G(V,fW)=fG(V,W),\nn\\
&&G(V,V)(q)\geq0,~~\mbox{for every }q\in M,\nn\\
&&G(V,V)(q)=0~~\mbox{implies that}~~V(q)=0.\nn
\eqn
A $(n,n)$-\ar\, is called a $n$-dimensional {\it almost-Riemannian structure}.
\edeff 
Let ${\cal S}=(M,\bD,G)$ be a  $(n,k)$-\ar. Reasoning as in  Remark~\ref{r-tensor}, we get that
$G$ is a tensor at the points $q$ where $\mathrm{dim}(\bD(q))=k$.
Although this is not necessarily the case everywhere on $M$, 
we can define,
for every $q\in M$, 
a quadratic form $\Gq$ on $\bD(q)$ through
\bqn
\Gq(v,v)=\inf\{G(V,V)(q)\mid V(q)=v, V\in\bD\}.
\eqnn

For every $q\in M$, it is possible to find a neighborhood $\Omega_q$ of $q$ and an {\it orthonormal frame} on $\Omega_q$, i.e., a 
set of  
$k$ vector fields $X_1,\dots, X_k\in \bD$ such that $G(X_i,X_j)=\delta_{i,j}$ on $\Omega_q$. One easily proves that 
orthonormal frames are local generators in $\Omega_q$.

If ${\cal S}$ is orientable then a volume form $\om$ can be chosen  
in such a way that $|\omega(X_1,\dots,X_k)|=1$ on every local orthonormal frame.

Let ${\cal S}=(M,\bD,G)$ be a  $(n,k)$-\ar. 
A curve $\g:[0,T]\to M$ is said to be {\it admissible} for ${\cal S}$ 
if  it is Lipschitz continuous and $\dot\g(t)\in\Delta_{\g(t)}$  for almost every $t\in[0,T]$.
Given an admissible 
curve $\g:[0,T]\to M$, the {\it length of $\g$} is  
\bqn
l(\g)= \int_{0}^{T} \sqrt{ \gg_{\gamma(t)}(\dot \g(t),\dot \g(t))}~dt.\eqnn
The {\it distance induced by ${\cal S}$ on $M$} is the function
\bqn
d(q_0,q_1)=\inf \{l(\g)\mid \g(0)=q_0,\g(T)=q_1, \g\ \mathrm{admissible}\}.
\eqnl{e-dipoi}
It is a standard fact that $l(\g)$ is invariant under reparameterization of the curve $\g$. 
Moreover, if an admissible curve $\g$ minimizes the so-called {\it energy functional} 
$
E(\g)=\int_0^T \gg_{\gamma(t)}(\dot \g(t),\dot \g(t))~dt
$
with $T$ fixed (and fixed initial and final point)
then $v=\sqrt{\gg_{\gamma(t)}(\dot \g(t),\dot \g(t))}$ is constant and 
$\g$ is also a minimizer of $l(\cdot)$. 
On the other side a minimizer $\g$ of $l(\cdot)$ such that  $v$ is constant is a minimizer of $E(\cdot)$ with $T=l(\g)/v$.

A {\it geodesic} for  ${\cal S}$  is a 
curve $\g:[0,T]\to M$ such that 
for every sufficiently small interval 
$[t_1,t_2]\subset [0,T]$, $\g|_{[t_1,t_2]}$ is a minimizer of $E(\cdot)$. 
A geodesic for which $\gg_{\gamma(t)}(\dot \g(t),\dot \g(t))$ is (constantly) 
equal to one is said to be parameterized by arclength.

The finiteness and the continuity of $d(\cdot,\cdot)$ with respect 
to the topology of $M$ are guaranteed by  the Lie bracket generating 
assumption on the \ar.  
The distance $d(\cdot,\cdot)$ gives to $M$ the 
structure of metric space.
The local existence of minimizing geodesics 
is a standard consequence of Filippov Theorem 
(see for instance \cite{agra-book}). 
When $M$ is compact 
any two points of $M$ are connected by a minimizing geodesic.

A convenient way to deal with a \ar\  is to identify it with an atlas of local orthonormal frames. In 
the case of an orientable 
\ar, one can impose that all orthonormal frames are coherently oriented. In this way one is led to the following 
equivalent definition.

\bdeff 
\llabel{d-ars}
Let $M$ be a $n$-dimensional smooth manifold, fix $k\in\N$, and consider a family
\bqn
{\cal S}=\{(\Omega^\mu,X_1^\mu,\dots,X_k^\mu)\}_{\mu\in I}, 
\eqnn
where $\{\Omega^\mu\}_{\mu\in I}$ is an
open covering of $M$ 
and, for every 
$\mu\in I$, $\{X_1^\mu,\dots,X_k^\mu\}$ is a family of 
smooth vector fields defined on $M$, whose restriction to 
$\Omega^\mu$ 
satisfies the Lie bracket generating condition. 
We  assume moreover that for every $\mu\in I$ and every open nonempty subset $\Omega$ of $\Omega^\mu$, 
the submodule of $\mathrm{Vec}(\Om)$ generated by $X_1^\mu,\dots,X_k^\mu$ on $\Omega$ cannot be 
generated by less than $k$ vector fields.

We say that ${\cal S}$ is a $(n,k)$-\ar\ if, 
for every $\mu,\nu\in I$ and for every $q\in \Omega^\mu\cap\Omega^\nu$, 
there exists an orthogonal matrix $R^{\mu,\nu}(q)=(R^{\mu,\nu}_{i,j}(q))\in 
O(k)$ such that 
\bqn
X^\mu_i(q)=\sum_{j=1}^kR^{\mu,\nu}_{i,j}(q)X^\nu_j(q).
\eqnl{e-theta}
We say that two \ar s ${\cal S}_1$ and ${\cal 
S}_2$ on $M$ are {\it equivalent} if ${\cal S}_1\cup {\cal S}_2$ is a
\ar. 
Given an open subset $\Omega$ of $M$ and a set of $k$ 
vector fields $(X_1,\dots,X_k)$, we say that  
$(\Omega,X_1,\dots,X_k)$ is {\it compatible 
with ${\cal S}$} if ${\cal S}\cup\{(\Omega,X_1,\dots,X_k)\}$ is equivalent to ${\cal 
S}$.

If ${\cal S}$ is equivalent to a \ar\ of 
the form $\{(M,X_1,\dots,X_k)\}$, i.e., for which the cardinality of $I$ is equal to one, we 
say that 
${\cal S}$ 
is \underline{trivializable}. 

If ${\cal S}$ admits an equivalent \ar\ such that each $R^{\mu,\nu}(q)$ belongs to $SO(k)$, 
we say that ${\cal S}$ is orientable.
\edeff

Notice that $R^{\mu,\nu}(q)$ is uniquely defined by equation \r{e-theta} and, 
moreover, is  smooth as a function of $q$. 
In the following, when dealing with an orientable \ar, we always assume that the  atlas of local orthonormal frames
is {\it positive oriented}, i.e., such that each $R^{\mu,\nu}$ belongs to $SO(k)$. For such an atlas, 
a volume form $\omega$ can be 
chosen such that $\omega(X^\mu_1,\ldots,X^\mu_k)=1$ on $\Omega^\mu$, for every $\mu\in I$.

In terms of Definition \ref{d-ars-1}, 
$\bD$ is the module that is locally (in $\Omega^\mu$)  generated  by  $X_1^\mu,\dots,X_k^\mu$. Moreover
\bqn
\bD(q)&=&\mathrm{span}\{X_1^\mu(q),\dots,X_k^\mu(q)\},\nn\\
\gg_q(v,v)&=&\inf\left\{\sum_{i=1}^k \al_i^2\mid v=\sum_{i=1}^k \al_i X^\mu_i(q)\right\},\nn 
\eqn
for every $q\in\Omega^\mu$ and $v\in\bD(q)$.

\bdeff
A property $(P)$ defined for $(n,k)$-\ar s on $M$ is said to be {\it generic}
if there exists an open and dense subset  ${\cal O}$ of   $(\mathrm{Vec}(M),{\cal W})^k$ such 
that $(P)$ holds for every $(n,k)$-\ar\ admitting an atlas of local orthonormal 
frames whose elements belong to 
${\cal O}$.
\edeff

Given a $(n,k)$-\ar\ ${\cal S}$, the problem of finding a curve minimizing the energy between two fixed points  $q_0,q_1\in M$ is 
naturally formulated as the optimal control problem
\bqn
&&\dot q=\sum_{i=1}^k u_i X^\mu_i(q)\,,~~~u_i\in\R\,,~~~\mu\in I(q)=\{\mu\in I\mid q\in\Omega^\mu\},\llabel{e-dyn}\\
&&\int_0^T 
\sum_{i=1}^k u_i^2(t)~dt\to \min,~~q(0)=q_0,~~~q(T)=q_1.\llabel{e-cost}
\eqn
Here $\mu,u_1,\dots,u_k$ are seen as controls and $T$ is fixed. 
It is a standard fact that this optimal control problem  is equivalent to the minimum time problem with 
controls $u_1,\ldots, u_k$ satisfying $u_1^2+\cdots+u_k^2\leq 1$.
	
Notice that if  the \ar\ is trivializable, then the role 
of 
$\mu$ is empty and  \r{e-dyn}, \r{e-cost} can be rewritten as a classical sub-Riemannian control problem 
\bqn
\dot q=\sum_{i=1}^k u_i X_i(q)\,,~~~u_i\in\R\,,
~~~\int_0^T 
\sum_{i=1}^k u_i^2(t)~dt\to\min,~~q(0)=q_0,~~~q(T)=q_1.
\eqnn

\section{Two-dimensional almost-Riemannian structures}\label{s-dim2}
Henceforth the paper is focused on the special case of $(2,2)$-\ar s, i.e., two-dimensional almost-Riemannian structures (2-ARSs for 
short). In this case a local orthonormal frame on $\Omega^\mu$ is a pair of Lie bracket generating vector fields $(X^\mu,Y^\mu)$. 

Given a 2-ARS ${\cal S}$, we
call {\it singular locus} the set
$\Zz\subset M$ 
of points $q$ at which the dimension of $\Delta(q)$ is equal to one.
Denote by $g$ the restriction of the quadratic form $\gg$ on $M\setminus\Zz$. By construction $g$ is a 
Riemannian metric satisfying 
\bqn
g(X^\mu(q),X^\mu(q))=1,~~
g(X^\mu(q),Y^\mu(q))=0,~~
g(Y^\mu(q),Y^\mu(q))=1,
\eqnn
for every $\mu$ in $I$ and every $q\in \Omega^\mu\setminus\Zz$.      
Denote moreover by $dA$ the Riemannian density associated with $(M\setminus\Zz,g)$, which coincides with
$|dX^\mu\wedge dY^\mu|$ on $\Om^\mu\setminus\Zz$, for every $\mu\in I$.

Finally, one can define on $M\setminus\Zz$ the Gaussian curvature $K$ associated with $g$, 
which  is easily expressed in each open set $\Omega^\mu\setminus\Zz$ 
through the formula (see for instance 
\cite{agra-book}, equation (24.6))
\bqn
K=-(\al^\mu)^2-(\beta^\mu)^2 +X^\mu \beta^\mu -Y^\mu\al^\mu, 
\eqnn
where $\al^\mu,\beta^\mu:\Omega^\mu\setminus\Zz\to\R$ are (uniquely) defined by
\bqn
[X^\mu,Y^\mu]=\al^\mu X^\mu+\beta^\mu Y^\mu,
\eqnn
and $X^\mu\beta^\mu$ (respectively, $Y^\mu\al^\mu$) denotes the Lie derivative of 
$\beta^\mu$ with respect to $X^\mu$ (respectively, of $\al^\mu$ with respect to 
$Y^\mu$).

\subsection{Minimizers, cut and conjugate loci}
\label{s-minimizers}
A natural tool to look for geodesics in almost-Riemannian geometry 
is to apply the necessary condition for optimality given by the 
Pontryagin Maximum Principle (see \cite{pont-libro}). As a result we 
obtain the following proposition. In view of later applications in the 
paper, we consider as initial and final conditions not only points, 
but submanifolds as well.
\bp
\llabel{p-pmp}
Define on $T^\ast M$ the Hamiltonian 
\bqn
H(\lam,q)=\frac12(\ll \lam,X^\mu(q)\rr^2+\ll \lam,Y^\mu(q)\rr^2),~~~~
\mbox{$q\in\Omega^\mu, ~~\lam\in T^\ast_q M$}.
\eqnn
(Notice that $H$ is well defined on the whole $T^\ast M$, thanks to 
\r{e-theta}.)
Consider the minimization problem 
\bqn
\dot q\in\bD(q),~~\int_0^T \gg_{q(t)}(\dot q(t),\dot q(t))~dt\to \min, 
~~q(0)=M_{\mathrm{in}},~~~q(T)=M_{\mathrm{fin}},
\eqnl{frog}
where $M_{\mathrm{in}}$ and $M_{\mathrm{fin}}$ are two submanifolds of $M$ and the final time $T>0$ is 
fixed.
Then every solution of \r{frog} is the projection on $M$ of a trajectory 
$(\lam(t),q(t))$ of the Hamiltonian system associated with $H$ 
satisfying 
$\lam(0)\perp T_{q(0)}M_{\mathrm{in}}$, $\lam(T)\perp 
T_{q(T)}M_{\mathrm{fin}}$, and $H(\lam(t),q(t))\neq0$. 
\ep
\brem
The simple form of the statement above follows from the absence of 
abnormal minimizers, which follows from the Lie bracket generating 
assumption. As a consequence a curve is a geodesic if and only if it the  projection 
of a normal extremal. 
\erem
\brem
Notice that $H$ is constant along any given  solution of the Hamiltonian system. Moreover, $H=1/2$ if and only if 
$q(.)$ is parameterized by arclength.
\erem
Fix $q\in M$. For every $\lam\in T_{q}^\ast M$ satisfying 
\bqn 
H(\lam, q)=1/2
\eqnl{1/2}
and 
every $t>0$ define $E(\lam,t)$ as the projection on $M$ of the solution, evaluated at time $t$, of 
the Hamiltonian system 
associated with $H$,  with initial condition $\lam(0)=\lam$ and $q(0)=q$.
Notice that if $q\notin\Zz$ then condition
\r{1/2} defines
an ellipse in  $T_{q}^\ast M$;  otherwise it 
identifies
a pair of parallel straight lines.

\bdeff
The \underline{conjugate locus from $q$} is the set of critical values of the map  $E(\lam,t)$.
For every $\bar\lam$ such that \r{1/2} holds, let $t(\bar\lam)$ be the first positive time, if it exists, 
for which the map $(\lambda,t)\mt E(\lam,t)$ is singular 
at $(\bar\lam,t(\bar\lam))$. The \underline{first conjugate locus} from 
$q$ is the set $\{E(\bar\lam,t(\bar\lam))\mid t(\bar\lam) \mbox{ exists}\}$.
The \underline{cut locus} from $q$ is the set of points reached optimally by more than one geodesic, 
i.e., the set
\bqn
\{q'\in M \mid \exists~\lam_1,\lam_2,t\mbox{ such 
that }q'=E(\lam_1,t)=E(\lam_2,t),~\lam_1\neq\lam_2,
\mbox{ and } E(\lam_1,\cdot), 
E(\lam_2,\cdot)\mbox{ are optimal in }[0,t]\}.  
\eqnn
\edeff
\brem
It is a standard fact that for every $\bar\lam$ satisfying \r{1/2}, the set 
$T(\bar\lam)=\{\bar t>0\mid\mathrm{~the~~map~~}(\lambda,t)\mt E(\lam,t)\mathrm{ ~~is~~  singular~~at~~}(\bar\lambda,\bar t) \}$ is a 
discrete set (see for instance \cite{agra-book}).
\erem

\ppotR{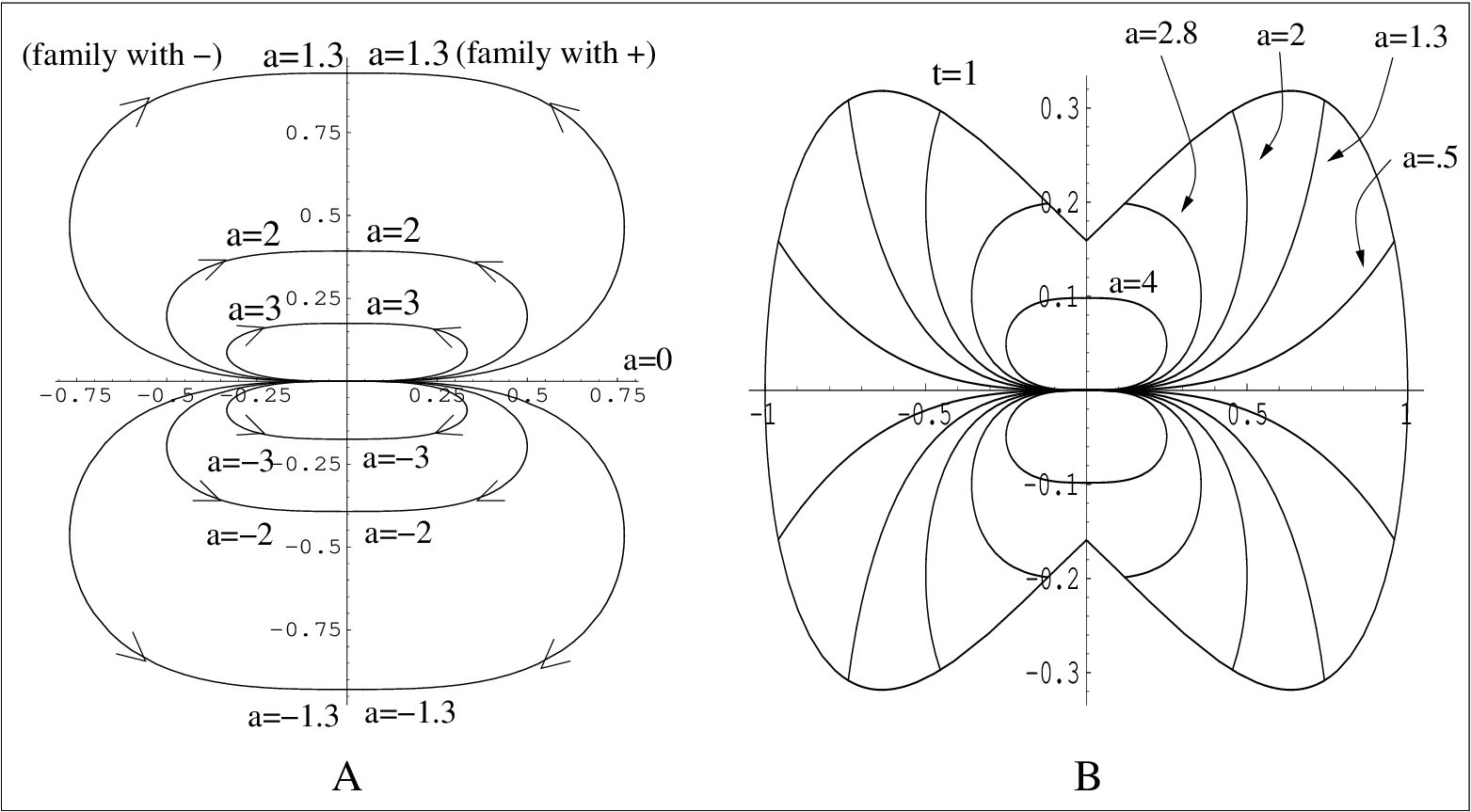}{Geodesics and minimum time front (for $t=1$) for 
the 
Grushin metric}{14}
\subsection{An example: the Grushin almost-Riemannian structure}
\llabel{ss-grushin}
Consider again the Grushin distribution $X(x,y)=(1,0)$, $Y(x,y)=(0,x)$ on
the plane $\R^2$. If we consider $X$ and $Y$ as an orthonormal frame, we get 
an almost-Riemannian structure.

As already remarked, the singular locus coincides with the $y$-axis. 
Therefore, every trajectory crossing the $y$-axis does it horizontally. 
The Riemannian metric $g$ associated 
with the Grushin metric on $\R^2\setminus\{(x,y)\in\R^2\mid x\neq0\}$
explodes when one is 
approaching the $y$-axis,
\bqn
g=dx^2+\frac{1}{x^2} dy^2.
\eqnn
Also the curvature and the Riemannian density explode while approaching the  $y$-axis,
\bqn
K=-\frac2{x^2}, ~~dA=\frac{1}{|x|}dx\,dy. 
\eqnn

According to Proposition \ref{p-pmp}, geodesics 
are the projection on the $(x,y)$-plane of the solutions of the 
Hamiltonian system corresponding to
$$
H=\frac12(\lam_x^2+\lam_y^2 x^2).
$$
Fixing the initial condition  $x(0)=0$, 
$y(0)=0$, the normalization $H=1/2$ implies that  
$\lam_x(0)=\pm1$. Taking $a=\lam_y(0)\in\R$, the geodesics starting from the origin 
are
\bqn
&&\mbox{for $a=0$}~~~~~
\left\{\ba{lll}
x_0(t)&=&\pm t\\
y_0(t)&=&0,\\
\ea\right.\nn\\
&&\mbox{for $a\neq0$}~~~~~\left\{\ba{lll}
x_a(t)&=&\pm\frac1{a}\sin(at)\\
y_a(t)&=&\frac{1}{2a} t-\frac{1}{4a^2}\sin(2at).\\
\ea\right.
\eqnn

Due to the symmetries of the problem, one can easily check that the 
time at which a geodesic $(x_a(t),y_a(t))$ loses optimality is  
$\bar t=\pi/|a|$, for $a\neq0$, and that $(x_a(\bar t), y_a(\bar t))$ belongs to the 
$y$-axis. The geodesics corresponding to $a=0$ are optimal for every positive time. 
As a consequence the cut locus from the origin is the set 
$\{(0,\al)\mid\al\in\R\setminus\{0\}\}$.

In Figure \ref{f-grusinGROSSA}A geodesics for some values of 
$a$ 
are portrayed, while Figure \ref{f-grusinGROSSA}B illustrates the set of 
points reached in time $t=1$. Notice that this set has a non-smooth boundary.  
In contrast with what would happen in Riemannian geometry, this is the case for every 
positive time, as it happens in constant-rank sub-Riemannian geometry. 
However, this is a consequence of the fact that the initial 
condition belongs to $\Zz$.

One can check that, even if the curvature is always negative 
where it is defined, a geodesic $(x_a(t),y_a(t))$, $a\neq0$, 
has its first conjugate point at time $\tau/|a|$, where $\tau\sim4.49$ is the 
first positive root of the equation $\tan(\tau)=\tau$. As a consequence 
the first conjugate locus is the parabola 
\bqn
y=\frac{x^2}{2}\left(\frac{1}{\cos\tau\sin\tau}-\frac{1}{\tau}\right).
\eqnn

One could ask whether the presence of conjugate points is the consequence 
of the particular initial point on the set $\Zz$. In fact this is not the 
case. Consider  as initial condition the point $x(0)=-1$, $y(0)=0$. 
Define, for every $a\in[0,1]$,
\bqn
&&x^+(t,a)=\frac{-\left( a\,\cos (a\,t) \right)  + 
    {\sqrt{1 - a^2}}\,\sin (a\,t)}{a},\nn\\
&&x^-(t,a)=\frac{-\left( a\,\cos (a\,t) \right)  - 
    {\sqrt{1 - a^2}}\,\sin (a\,t)}{a},\nn\\
&&y^+(t,a)=\frac{-4\,a\,{\sqrt{1 - a^2}} + 2\,a\,t + 
    4\,a\,{\sqrt{1 - a^2}}\,{\cos (a\,t)}^2 - 
    \sin (2\,a\,t) + 2\,a^2\,\sin (2\,a\,t)}{4\,a^2},\nn\\
&&y^-(t,a)=\frac{4\,a\,{\sqrt{1 - a^2}} + 2\,a\,t - 
    4\,a\,{\sqrt{1 - a^2}}\,{\cos (a\,t)}^2 - 
    \sin (2\,a\,t) + 2\,a^2\,\sin (2\,a\,t)}{4\,a^2}.\nn
\eqn
Then every geodesic  from the point $(-1,0)$ belongs to one of the four families
\bqn
&&G1:~~~(x^+(t,a),y^+(t,a)),\nn\\
&&G2:~~~(x^-(t,a),y^-(t,a)),\nn\\
&&G3:~~~(x^+(t,a),-y^+(t,a)),\nn\\
&&G4:~~~(x^-(t,a),-y^+(t,a)).\nn
\eqn
The geodesics in G1 and G2 lie in the half plane 
$\{y\geq0\}$, while those in G3 and G4 lie in $\{y\leq0\}$.

Let us describe the cut locus from $(-1,0)$.  Consider first a geodesics $(x^+(t,a),y^+(t,a))$
belonging to the family G1.
One can check that $\pi/a$ is the first positive
time  at which $(x^+(t,a),y^+(t,a))$
intersects another geodesic, namely  $(x^-(t,a),y^-(t,a))$, which belongs to G2.
The situation is similar for the 
families G3 and G4. As a consequence the cut locus 
from $(-1,0)$ is the set 
$\{(1,\al)\mid\al\in[\pi/2,\infty)\cup(-\infty,-\pi/2]\}$. 

As 
above,
 one can also check that every geodesic 
(except those corresponding to $a=0$) has a conjugate time (see 
Figure \ref{f-cong-cut}).
In particular $\pi$ is a conjugate time for the geodesics 
corresponding to $a=1$. 
Notice that conjugate 
points appear on geodesics which have already crossed $\Zz$.
(Before crossing $\Zz$ a geodesic is Riemannian and 
lies in a Riemannian space with negative Gaussian curvature.)

\ppotR{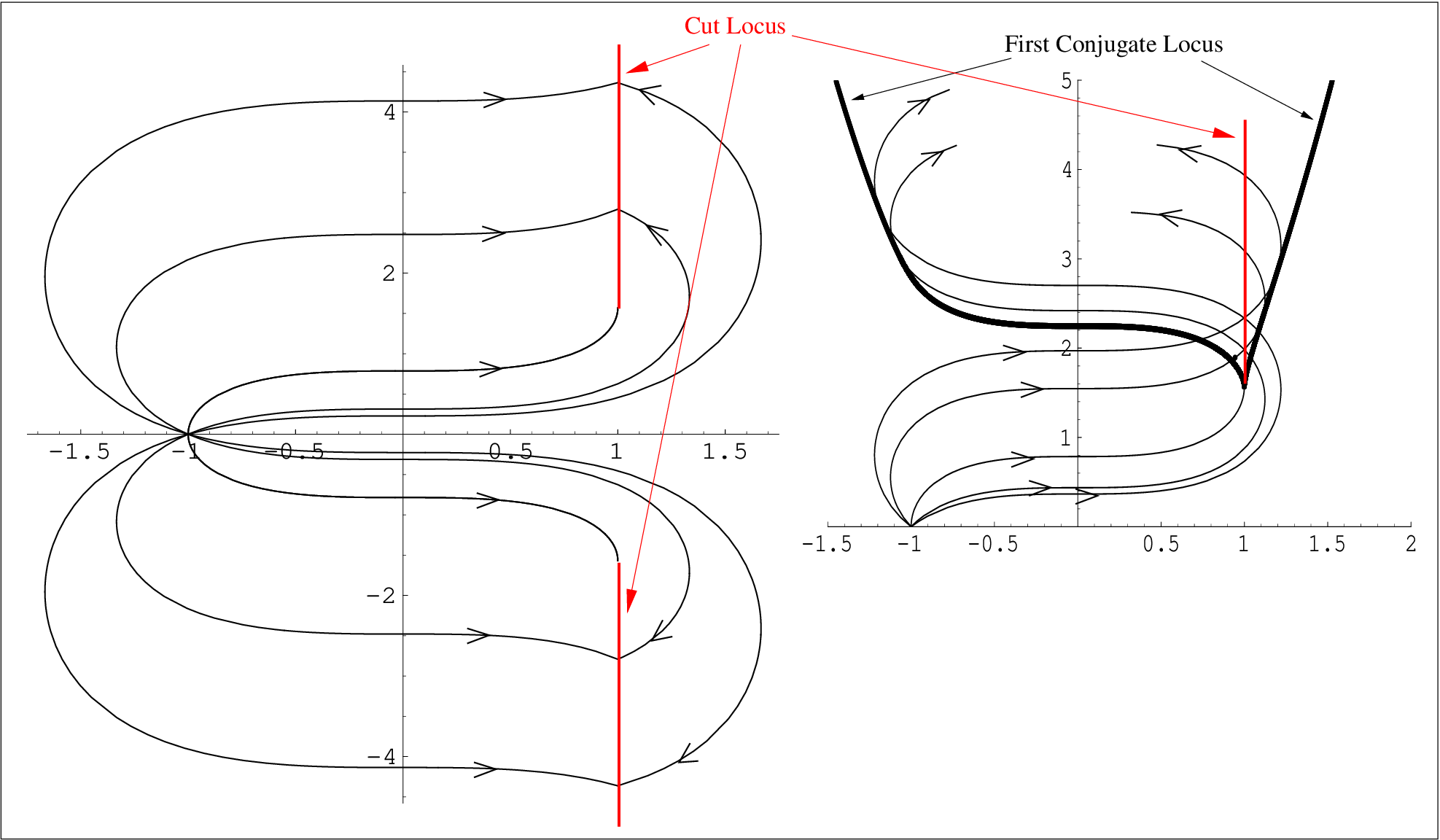}{}{17}

\section{Normal forms for generic 2-ARSs}
\label{s-generic}

The following proposition is a standard corollary of the transversality theorem. It formulates generic properties of a 2-ARS in terms of the 
flag of the distribution $\Delta$ (see equation \r{flag}).

\bp \label{p-generic}
Let $M$ be a two-dimensional smooth manifold. Generically, a 2-ARS ${\cal 
S}=\{(\Omega^\mu,X^\mu,Y^\mu)\}_{\mu\in I}$ on $M$ satisfies the following 
properties: {\bf (i)} $\Zz$ is an 
embedded one-dimensional smooth  
submanifold of 
$M$; 
{\bf (ii)} The points $q\in M$ at which $\bD_2(q)$ is 
one-dimensional are isolated;  
{\bf (iii)}  $\bD_3(q)=T_qM$ for every $q\in M$.
\ep
\brem
\label{r-asub}
Notice that properties 
{\bf (i)},
{\bf (ii)}, and 
{\bf (iii)} are actually generic for every $(2,2)$-\rvd, since they do not involve the metric 
structure. 
\erem
As a consequence of Proposition~\ref{p-generic}, one can classify the 
local normal forms of a generic 2-ARS.
\bt
\label{t-normal}
Generically for a  2-ARS ${\cal S}$, for every point 
$q\in M$ there exist a neighborhood $U$ of $q$ and a pair of vector 
fields $(X,Y)$ on $M$ such that $(U,X,Y)$ is compatible with ${\cal S}$ and, up 
to a smooth change of coordinates defined on $U$, $q=(0,0)$ and $(X,Y)$ 
has one of the 
forms
\bqn
\mathrm{(F1)}&& ~~X(x,y)=(1,0),~~~Y(x,y)=(0,e^{\phi(x,y)}), \nn  \\
\mathrm{(F2)}&& ~~X(x,y)=(1,0),~~~Y(x,y)=(0,x e^{\phi(x,y)}),\nn   \\
\mathrm{(F3)}&& ~~X(x,y)=(1,0),~~~Y(x,y)=(0,(y -x^2 
\psi(x))e^{\phi(x,y)}), \nn
\eqn
where $\phi$ and $\psi$ are smooth real-valued functions such that
$\phi(0,y)=0$ and  $\psi(0)\neq0$.
\et
Before proving Theorem \ref{t-normal} let us show the following lemma.
\bl
\label{l-w}
Let ${\cal S}$ be a 2-ARS and let
$W$ be a smooth embedded one-dimensional submanifold of $M$. Assume that $W$ is 
transversal to 
the distribution $\bD$, i.e., such that 
$\bD(q)+T_qW=T_qM$
 for every $q\in W$.
Then, for every $q\in W$ there exists an open  
neighborhood $U$ of $q$ 
such that for every $\eps>0$ the set 
\bqn
\{q'\in U\mid d(q',W)=\eps\},
\eqnn
is a smooth  embedded one-dimensional submanifold of $U$. Moreover, there 
exists a pair of vector fields $(X,Y)$ such that $(U,X,Y)$ is compatible with ${\cal S}$ and,  up to a 
smooth change of coordinates defined on $U$, $q=(0,0)$ and 
$W, X$, $Y$  
have the form
\bqn
W&=&\{(0,h)\mid h\in\R\},\nn\\
X(x,y)&=&(1,0),\nn\\
Y(x,y)&=&(0,f(x,y))\nn
\eqn
where $f(x,y)$ is a smooth function defined on $U$.
\el
\proof
Consider a smooth regular parametrization $\al\mapsto w(\al)$ of 
$W$. Let $\al\mapsto \lam_0(\al)\in T^\ast_{w(\al)}M$ be a smooth map 
satisfying  $H(\lam_0(\al),w(\al))=1/2$ and 
$\lam_0(\al)\perp T_{w(\al)}W$. 

Let $E(t,\al)$ be the solution at time $t$ of the Hamiltonian system given 
by the Pontryagin Maximum Principle with initial condition 
$(q(0),\lam(0))=(w(\al),\lam_0(\al))$ (see Proposition \ref{p-pmp}).
Fix $q\in W$ and define $\bar\al$ by $q=w(\bar\al)$.
In order to prove that $E(t,\al)$ is a local diffeomorphism around the 
point $(0,\bar\al)$, let us show that the two vectors
\bqn
v_1=\frac{\partial E}{\partial \al}(0,\bar\al)\mbox{~~and }
v_2=\frac{\partial E}{\partial t}(0,\bar\al)
\eqnn
are not parallel. On one hand, since $v_1$ is equal to $\frac{dw}{d\al}(\bar\al)$, then it   
spans $T_qW$. 
On the other hand, 
being $H$ quadratic in 
$\lam$, 
\bqn
\ll \lam_0(\bar\al),v_2\rr= 
\ll \lam_0(\bar\al),\frac{\partial H}{\partial \lam}(\lam_0(\bar\al),q)\rr=
2 H(\lam_0(\bar\al),q)=1.
\eqnn
Thus $v_2$ does not belong to the orthogonal to $\lam_0(\bar \al)$, that 
is, to $T_qW$.  

Therefore for a small enough neighborhood $U$ of $q$ 
we have that the set
$\{q'\in U\mid d(q',W)=\eps\}$ is given by the 
intersection of $U$ with the images of 
$E(\eps,\cdot)$
and $E(-\eps,\cdot)$. This proves the first part of the statement. 
To prove the second part, let us take $(t,\al)$ as a system of coordinates on $U$ 
and define the vector field $X$ by
\bqn
X(t,\al)=\frac{\partial E(t,\al)}{\partial t}.   
\eqnn
Notice that, by construction, for every $q'\in U$ the vector  $X(q')$ belongs to $\bD(q')$ and  
$\gg_{q'}(X(q'),X(q'))=1$. In the coordinates $(t,\al)$ we have 
$X=(1,0)$. Let $Y$ be a vector field on $U$  such that $(X,Y)$ is 
compatible with ${\cal S}$. 
We are left to prove that the first component of $Y$ is identically equal to zero. 
Indeed, were  
this not the case, the norm of $X$ would not be equal to one.
\quadp\\\\
{\bf Proof of Theorem \ref{t-normal}.}
Let us start from the case in which $\Delta(q)=T_q M$. Let $W$ be any 
one-dimensional submanifold passing through $q$. Lemma \ref{l-w} provides us 
with a possible choice of orthonormal frame $X=(1,0)$, $Y=(0,f(x,y))$ in a 
neighborhood $U$ of $q$.
Since, without loss of generality, 
$X$ and $Y$ are everywhere linearly independent in $U$, then 
$f(x,y)\neq0$ for every $(x,y)\in U$. 
By applying a smooth coordinate transformation of the type $x\to x$, $y\to\nu(y)$ we get the 
new expressions $X=(1,0)$, $Y=(0,\nu'(y)f(x,y))$. A normal form of 
type (F1) is obtained by choosing $\nu$ in such a 
way that $\nu'(y)f(0,y)=1$.

Let now $q\in\Zz$ and assume that $\bD_2(q)=T_q M$. Assume, moreover, that 
the generic condition {\bf (i)}
holds true. One can easily check  that 
$\bD(q)$ is transversal to the submanifold $\Zz$ at $q$. Hence we can apply 
Lemma \ref{l-w} with $W=\Zz$. As a result we obtain  
a possible choice of orthonormal frame 
$X=(1,0)$, $Y=(0,f(x,y))$ in a
neighborhood $U$ of $q$. Since $X$ and $Y$ are linearly dependent on $\Zz$, 
which is identified with the $y$-axis,
then $f(0,y)=0$. 
The condition  $\bD_2(q)=T_q M$ implies that, by taking $U$ small enough,  
$\partial_x f(0,y)\neq0$. Hence $f$ admits a representation of the 
type $f(x,y)=x e^{\phi(x,y)}$, with $\phi$ smooth. 
Again, a change of coordinates $x\to x$, $y\to\nu(y)$ 
can be used in order to ensure that $\phi(0,y)=0$. 
The normal form 
(F2) is obtained.

Let  now $q\in\Zz$ be such that  $\bD_2(q)=\bD(q)$. Assume that
the 
generic conditions  {\bf (i)}, {\bf (ii)}, {\bf (iii)} are fulfilled. Let $W$ be any 
one-dimensional submanifold passing through $q$ and being 
transversal
to 
$\Zz$. Using Lemma  \ref{l-w} we can chose $X=(1,0)$, $Y=(0,f(x,y))$ with
$f$ satisfying, by assumption, 
$\partial_x f(0,0)=0$, $\partial^2_x f(0,0)\neq0$. Let us identify $\Zz$ 
with the graph of a smooth function $y=\Gamma(x)$. Then 
$f(x,y)$ can be written in the form $(y -\Gamma(x))e^{\phi(x,y)}$ with 
$\phi$ smooth. As above, without loss of generality $\phi(0,y)=0$. 
The conditions on $f$ at $(0,0)$, moreover, justify
the representation
$\Gamma(x)=x^2\psi(x)$, with $\psi$ smooth.
\quadp

\brem
Because of Remark \ref{r-asub}, 
for a generic distribution $\bD$, every 2-ARS having $\bD$ as 
corresponding distribution
can be locally represented by one of the normal forms 
(F1), (F2), (F3).
\erem

\bdeff
Let ${\cal S}$ be a 2-ARS and assume that the generic conditions
{\bf (i)}, {\bf (ii)},  {\bf (iii)} of Proposition \ref{p-generic} hold true. 
A point $q\in M$ is said to be an \underline{ordinary point} 
if $\bD(q)=T_q M$, hence, if 
${\cal S}$ is locally described by (F1). 
We call $q$ a  \underline{Grushin point} if $\bD(q)$ is one-dimensional and $\bD_2(q)=T_q M$, i.e., if  the local description (F2) applies. 
Finally, if $\Delta(q)=\Delta_2(q)$ is of dimension one and $\bD_3(q)=T_q M$
we say that $q$ is a \underline{tangency point} and ${\cal S}$ can be described near $q$ 
by the normal form 
(F3). \edeff
Let us take advantage of the common expressions of  the normal forms (F1), 
(F2), (F3), which are all of
the type $X(x,y)=(1,0)$, $Y(x,y)=(0,f(x,y))$, in order to investigate
the local
behavior of 
$g$, $K$, and $dA$.
\bl
\llabel{l-diagonal}
Let $X(x,y)=(1,0)$ and $Y(x,y)=(0,f(x,y))$ be two smooth vector fields on 
$\R^2$. 
Let $D=\{(x,y)\in \R^2~|~~f(x,y)\neq0 \}$ and $g$ be the Riemannian metric 
on $D$ having $(X,Y)$ as an orthonormal frame. Denote by $K$ the
 curvature of $g$ and by $dA$ the Riemannian density. We 
have
\bqn
g&=&dx^2+\frac1{f^2} dy^2,\nn\\
K&=&\frac{-2\,({\partial_x f})^2 + f\,\partial^2_x 
f}{{f}^2},\nn\\
dA&=&\frac1{|f|}\,dx\,dy.
\eqnn
\el
\subsection{An example of tangency point}
A simple example of tangency point can be observed in the 2-ARS defined by 
\bqn
X(x,y)=(1,0),~~Y(x,y)=(0,y-x^2).
\eqnn

The pair $(X,Y)$ 
appears in the classification of planar phase portraits of pairs of 
vector fields given by Davydov 
in \cite{davidov-book}. 
For this system one has
$g=dx^2+(y-x^2)^{-2} dy^2$, and
\bqn
K=\frac{-2\,\left( 3\,x^2 + y \right) }
  {{\left( x^2 - y \right) }^2}.
\eqnn

The graph of $K$ is illustrated in 
Figure 
\ref{curvature-tang}.
Notice that, in contrast with the behavior of the curvature in the 
Grushin plane (see Section \ref{ss-grushin}), in this case 
$\limsup_{q\to (0,0)} 
K(q)=+\infty$, while we still have $\liminf_{q\to (0,0)} K(q)=-\infty$. 
\ppotR{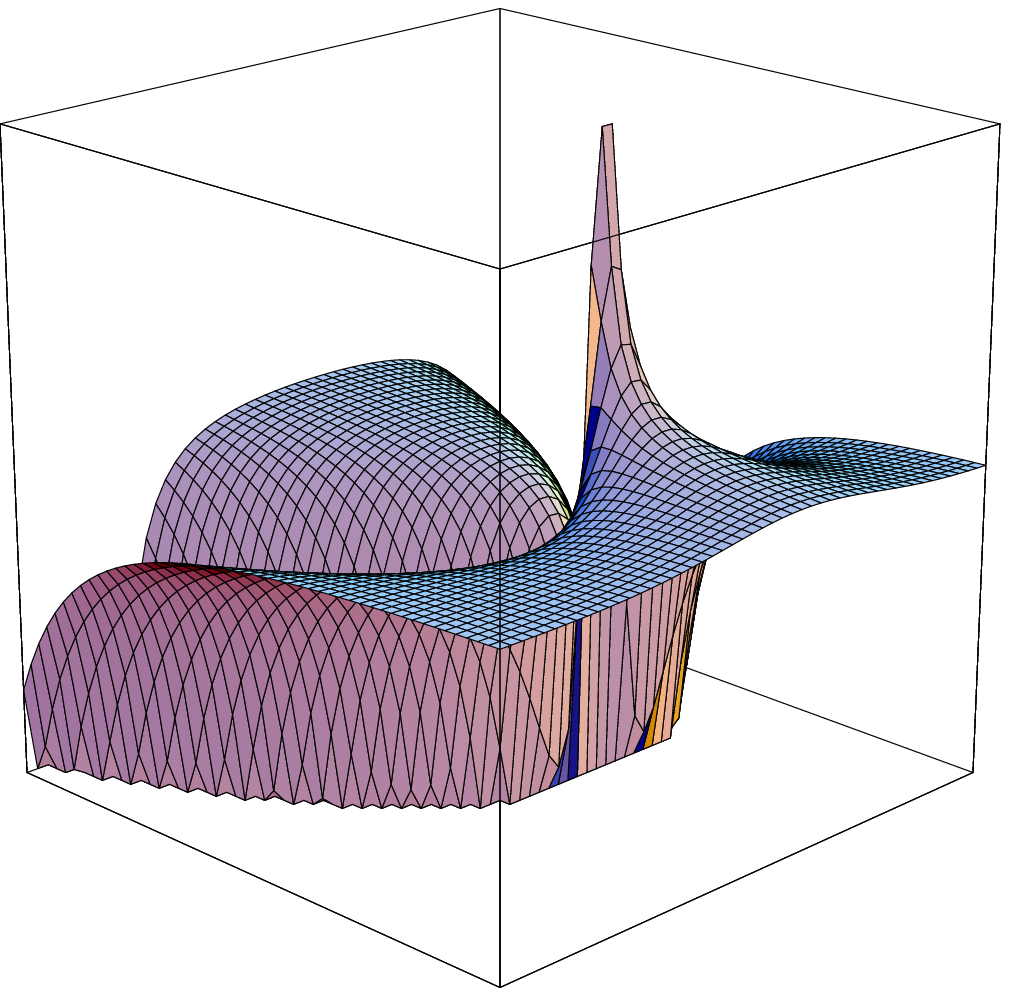}{}{10}

\section{The main result}
\label{s-main}
\subsection{Statement}

Let $M$ be an orientable two-dimensional manifold and let ${\cal S}$ be an orientable 2-ARS on $M$.
Chose a positive oriented atlas of orthonormal frames $\{(\Omega^\mu,X^\mu,Y^\mu)\}_{\mu\in I}$ 
and  denote by $\om$ the volume form  
on ${\cal S}$ such that $\om(X^\mu,Y^\nu)=1$ on $\Om^\mu$ for every $\mu\in I$.
As noticed in Remark~\ref{r-tensor}, $\om$ acts as a tensor on $M\setminus\Zz$. 
Define a two-form $dA_s$ on $M\setminus\Zz$ 
by the rule 
$dA_s(V(q),W(q))=\om(V,W)(q)$.
Notice that $dA_s=dX^\mu\wedge dY^\mu$ on $\Omega^\mu\setminus\Zz$ for every 
$\mu\in I$. 

Fix now an orientation $\Xi$ of $M$.
Recall that the choice of $\Xi$ determines uniquely a notion of integration on 
$M\setminus \Zz$ with respect to the form $dA_s$.  More precisely, given a  $dA$-integrable function $f$ on $\Om\subset M$, if for 
every $q\in\Omega$, $\Xi$ and $dA_s$ define the same orientation at $q$ (i.e. if $\Xi(q)=\al dA_s(q)$ with $\al>0$), then 
\bqn
\int_\Omega f ~dA_s=\int_{(\Omega,\Xi)}f~dA_s=\int_\Omega f ~|dA_s|=\int_\Omega f ~dA.
\eqnn
Let
\bqn
M^\pm=\{q\in \Om^\mu\setminus\Zz\mid \mu\in I,\pm\Xi(X^\mu,Y^\nu)(q)>0\}.
\eqnn
Then $\int_\Om f dA_s=\pm\int_\Om f dA$ if $\Om\subset M^\pm$.

For every $\eps>0$ let $M_\eps=\{q\in M\mid d(q,\Zz)>\eps\}$, where  $d(\cdot,\cdot)$ is 
the almost-Riemannian distance (see equation \r{e-dipoi}).  
We say that $K$ is {\it ${\cal S}$-integrable} if 
\bqn
\lim_{\eps\to0}\int_{M_\eps}K~dA_s
\eqnn
exists and is finite. In this  case we denote such limit by $\int K dA_s$.

\bt\llabel{t-1}
Let $M$ be a compact 
oriented 
two-dimensional manifold. For a generic oriented 2-ARS 
on $M$ such that no tangency point exists, $K$ is ${\cal S}$-integrable and
\bqn
\int K dA_s=2\pi
(\chi(M^+)-\chi(M^-)),
\eqnn
where  $\chi$ denotes the Euler characteristic.
\et
Theorem \ref{t-1} is proved in Section \ref{ss-prova}. 
For a generic trivializable 2-ARS without tangency points one can show, 
thanks to topological 
considerations (see Section \ref{s-topo}), 
that $\chi(M^+)=\chi(M^-)$. 
As a consequence, we derive the following 
result.
\bc\label{trivial}
Let $M$ be a compact oriented two-dimensional manifold. For a generic 
trivializable 2-ARS on $M$ without tangency points we have
\bqn
\int K dA_s=0.
\eqnn
\ec
\brem
In the results stated above, the hypothesis that there are not tangency points 
seems to be essential. Technically, the difficulty comes when one tries to 
integrate the Hamiltonian system given by the 
Pontryagin Maximum Principle applied to a system written in the normal form (F3).
However it is our hope to extend the Gauss-Bonnet formula even in presence of 
tangency points, using a more general approach.
\erem
It is anyway interesting to notice that the hypotheses of
Corollary \ref{trivial} are never empty, independently of $M$. Indeed:
\bl
\llabel{l-existence}
Every compact orientable 
two-dimensional manifold admits a trivializable 2-ARS satisfying the 
generic conditions of Proposition \ref{p-generic} and having  no 
tangency points. 
\el
The proof of Lemma \ref{l-existence} is given in Section \ref{s-sketch}.





\subsection{Proof of Theorem \ref{t-1}}
\llabel{ss-prova}
As a consequence of the compactness of $M$ and of Lemma~\ref{l-w} one easily gets:
\bl
\llabel{l-anelli} 
Let $M$ be compact and oriented. For a generic 2-ARS ${\cal S}$ on $M$ 
the set  
$\Zz$ is 
the union of finitely many curves diffeomorphic to  $S^1$.  
Moreover,
there exists $\eps_0>0$ such that, for every 
$0<\eps<\eps_0$, the set  
$M\setminus M_\eps$ is 
homeomorphic to $\Zz\times [0,1]$.  
Under the additional assumption that $M$ contains no tangency point, 
$\eps_0$ can be taken in such a way that 
$\partial M_{\eps}$ is smooth for every 
$0<\eps<\eps_0$. 
\el

Fix $M$ and a 2-ARS ${\cal S}$ as in the statement of Theorem \ref{t-1}.
Thus, ${\cal S}$ can be described, around each point 
of $\Zz$, by a normal form of type (F2).

Take $\eps_0$ as in the statement of Lemma \ref{l-anelli}. 
For every $\eps\in (0,\eps_0)$, let $M^\pm_\eps=M^\pm\cap M_\eps$. 
By 
definition of $dA_s$ and $M^\pm$,
\bqn
\int_{M_\eps}K dA_s=
\int_{M_\eps^+}K dA-
\int_{M_\eps^-}K dA.
\eqnn
The Gauss-Bonnet formula asserts that for every compact oriented 
Riemannian manifold $(N,g)$ with smooth 
boundary $\partial N$, we have
\bqn
\int_N K dA+\int_{\partial N} k_g ds=2\pi\chi(N),
\eqnn
where $K$ is the curvature of $(N,g)$, $dA$ is the 
Riemannian density, $k_g$ is the geodesic curvature 
of $\partial N$ 
(whose orientation is induced by  the one of $N$), and $ds$ is the 
length element. 

Applying the Gauss-Bonnet formula to the Riemannian manifolds $(M^+_\eps,g)$ 
and $(M^-_\eps,g)$ (whose boundary smoothness is guaranteed by Lemma 
\ref{l-anelli}), we have
\bqn
\int_{M_\eps}K dA_s=
2\pi(\chi(M^+_\eps)-\chi(M^-_\eps))-
\int_{\partial M^+_\eps}k_g ds+   
\int_{\partial M^-_\eps}k_g ds.   
\eqnn
Thanks again to 
Lemma~\ref{l-anelli}, $\chi(M^\pm_\eps)=\chi(M^\pm)$. We are left to prove that 
\bqn
\lim_{\eps\to0}\left(\int_{\partial M^+_\eps}k_g ds-   
\int_{\partial M^-_\eps}k_g ds\right)=0.   
\eqnl{e-limite}

Fix $q\in \Zz$ and a (F2)-type local system of coordinates $(x,y)$ in a neighborhood $U_q$ of $q$.  
We can assume that $U_q$ is given, in the coordinates $(x,y)$, 
 by a rectangle $[-a,a]\times [-b,b]$, $a,b>0$. 
Assume that $\eps<a$. 
Notice that $\Zz\cap U_q=\{0\}\times[-b,b]$ and   $\partial M_\eps\cap 
U_q=\{-\eps,\eps\}\times[-b,b]$. 
\ppotR{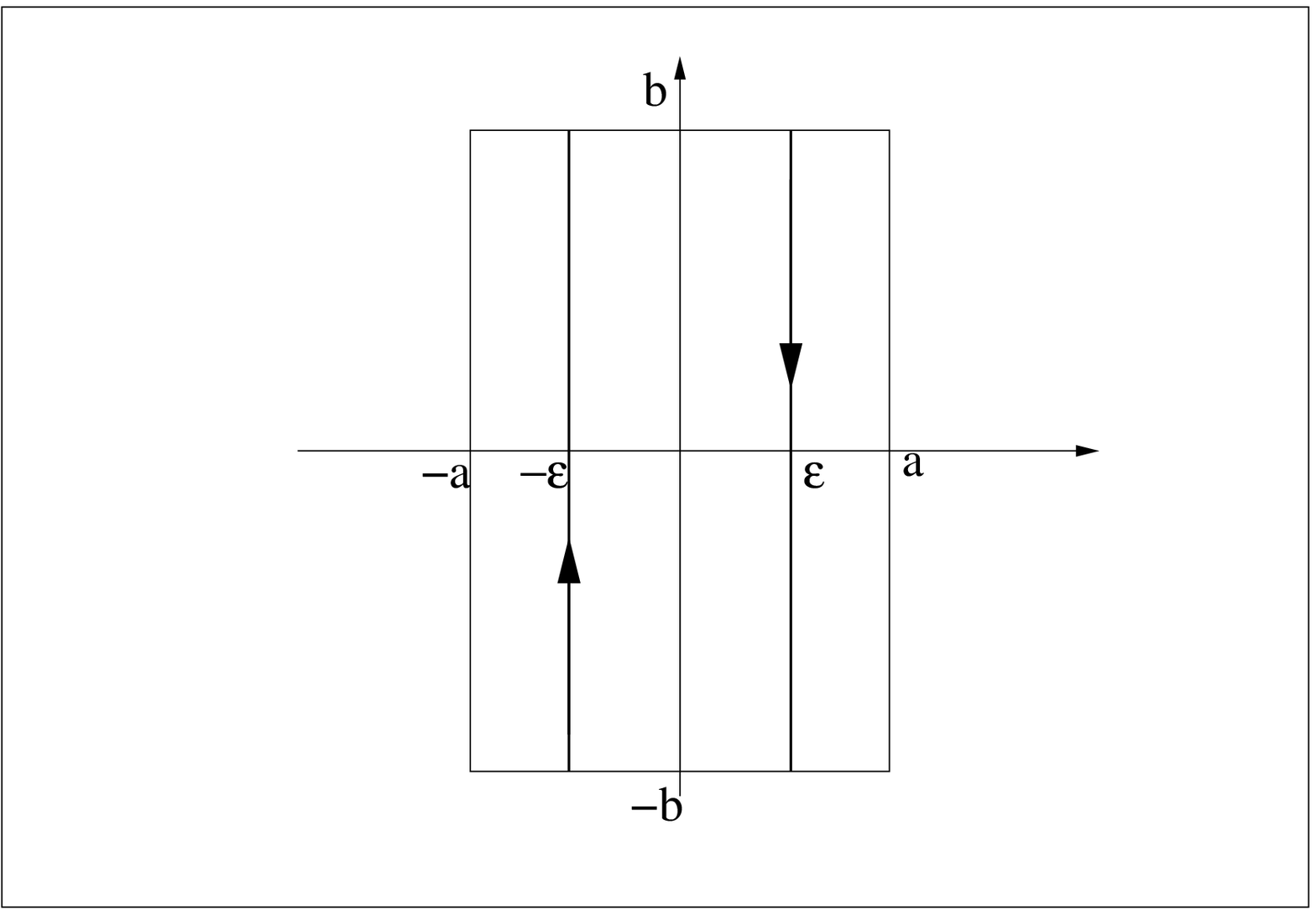}{}{8}
We are going to prove that 
\bqn
\int_{\partial M_\eps\cap U_q} k_g ~ds=O(\eps). 
\eqnl{e-integral-kg}
Then \r{e-limite} follows from the compactness of $\Zz$. (Indeed, $\{-\eps\}\times[-b,b]$ and $\{\eps\}\times[-b,b]$,  the horizontal edges of $\partial U_q$, are 
geodesics minimizing the length  from $\Zz$. Therefore, $\Zz$ can be covered by a finite number of neighborhoods 
of type $U_q$ whose pairwise intersections have empty interior.)

Without loss of generality, we can assume that $M^+\cap 
U_q=(0,a]\times[-b,b]$. Therefore,
$M^+_\eps$ induces on $\partial M^+_\eps=\{\eps\}\times[-b,b]$ a 
downwards orientation (see  
Figure \ref{f-rettangolo}). 
The curve $s\mt c(s)=(\eps,y(s))$ satisfying
\[\dot c(s)=-Y(c(s))\,,\ \ \ c(0)=(\eps,0)\,,\]
is an oriented parametrization by arclength of $\partial M^+_\eps$, making 
a constant angle with  $X$. 
Let $(\theta_1,\theta_2)$ be the dual basis to $(X,Y)$ on $U_q\cap M^+$,
i.e., $\theta_1=dx$ and $\theta_2=x^{-1} e^{-\phi(x,y)}dy$.
According to \cite[Corollary~3, p.~389, Vol.~III]{Spivak}, the geodesic 
curvature of $\partial M^+_\eps$ at $c(s)$ is equal to
$\lambda(\dot c(s))$,
where $\lambda\in\Lambda^1(U_q)$ is the unique 
one-form satisfying
\[ 
d\theta_1=\lambda\wedge \theta_2\,,\ \ \ d\theta_2=-\lambda 
\wedge\theta_1\,.
\]
A trivial computation shows that
\[\lambda=\partial_x(x^{-1}e^{-\phi(x,y)})dy\,.\]
Thus, 
\[k_g(c(s))=-\partial_x(x^{-1}e^{-\phi(c(s))})\,(dy(Y))(c(s))=
\frac1\eps+\partial_x\phi(\eps,y(s))\,.\]

Denote by $L_1$ and $L_2$ the lengths of, respectively, 
$\{\eps\}\times[0,b]$ and $\{\eps\}\times[-b,0]$. Then,
\begin{eqnarray*}
\int_{\partial M^+_\eps\cap U_q}k_g ds&=&\int_{-L_1}^{L_2}k_g(c(s))ds\\
&=&
\int_{-L_1}^{L_2} \lp \frac1\eps+\partial_x\phi(\eps,(s))\rp ds\\
&=&\int_{-b}^b \lp\frac1\eps+\partial_x\phi(\eps,y)\rp\frac1{\eps e^{\phi(\eps,y)}}dy\,,
\end{eqnarray*}
where the last equality is obtained taking $y=y(-s)$ as new variable of integration. 

We reason similarly on $\partial M^-_\eps\cap U_q$, on which $M^-_\eps$ 
induces the  
upwards orientation. An orthonormal 
frame on $M^-\cap U_q$,
oriented consistently with $M$, is given by 
$(X,-Y)$, whose dual basis is  
$(\theta_1,-\theta_2)$.  The same 
computations as above lead to 
\[\int_{\partial 
M^-_\eps\cap U_q}k_g ds=\int_{-b}^b \lp\frac1\eps-\partial_x\phi
(-\eps,y)\rp\frac1{\eps e^{\phi(-\eps,y)}}dy\,.\]
Define 
\be\label{d_F}
F(\eps,y)=(1+\eps\partial_x\phi(\eps,y))e^{-\phi(\eps,y)}.
\ee
Then 
\bqn
\int_{\partial 
M^+_\eps\cap U_q}k_g ds-
\int_{\partial 
M^-_\eps\cap U_q}k_g ds=\frac1{\eps^2}\int_{-b}^b (F(\eps,y)-F(-\eps,y))~dy.
\eqnn
By Taylor expansion with respect to $\eps$ we get
$$F(\eps,y)-F(-\eps,y)=2\partial_\eps F(0,y)\eps+O(\eps^3)=O(\eps^3)$$
where the last equality follows from the relation $\partial_\eps F(0,y)=0$ (see equation \r{d_F}).
Therefore,
\bqn
\int_{\partial 
M^+_\eps\cap U_q}k_g ds-
\int_{\partial 
M^-_\eps\cap U_q}k_g ds=O(\eps),
\eqnn
and \r{e-integral-kg} is proved. \quadp


\subsection{A counterexample in the non-generic case}
In this section we 
justify
the assumption that the 2-ARS is 
generic, by presenting an example of 2-ARS such that the conclusion 
of Theorem~\ref{t-1} does not hold, although $\Zz$ is smooth and $\bD(q)$ is transversal to $\Zz$ at every point $q$ of $\Zz$. 

Let $M$ be the two-dimensional torus $[-\pi,\pi]\times[-\pi,\pi]$ with the 
standard identifications, and consider the trivializable 2-ARS associated 
with the vector fields
\bqn
X(x,y)=(1,0),~~~
Y(x,y)=(0,1-\cos(x)).
\eqnn 
In this case $\Zz$ is the circle $\{0\}\times[-\pi,\pi]$ and one among 
$M^+$ and $M^-$ is empty (say $M^-$). Notice that the generic condition 
{\bf (ii)} is not verified since $[X,Y](q)=0$ at every $q\in \Zz$. By Lemma 
\ref{l-diagonal} we have
\bqn
dA&=&\frac1{1-\cos(x)}dx\,dy,\nn\\
K&=&\frac{\cos (x)-2 
       }{2\sin (\frac{x}{2})^2},\nn
\eqn   
 on $M\setminus\Zz$. Thus,
$\int K dA_s=\int_{M\setminus\Zz} K dA=-\infty$.


\subsection{Trivializable 2-ARSs}
\llabel{s-topo}

The aim of  this section is to characterize topologically 
trivializable 2-ARSs having no tangency point.

\bl
Let $M$ be orientable. 
For a generic trivializable 2-ARS on $M$ 
without tangency points
the Euler 
characteristics of $M^+$ and $M^-$ are equal.
\el
\proof
Let us consider on $M$ a notion of angle, induced 
by any fixed, globally defined, Riemannian metric $g_0$. 
For every 
 $\theta$ in $S^1$, denote by  $R_{\theta}:TM\rw TM$ the corresponding rotation 
of angle 
$\theta$.  

Since, by hypothesis, the map $\Zz\ni q\mapsto\bD(q)$
is a
one-dimensional distribution  
everywhere transversal to the smooth submanifold $\Zz$, 
then
we can define a 
smooth
function  $\phi:\Zz\to (0,\pi)$ such that
$R_{\phi(q)}(\Delta(q))=T_q\Zz$ for every $q\in\Zz$.
Let $\theta:M\times [0,1]\to S^1$ be a $\con^\infty$ map such that 
\bqn
&&\theta|_{M\times\{0\}}\equiv 0\nn\\
&&\theta|_{\Zz\times\{1\}}= \phi.\nn
\eqn
For every $t\in[0,1]$ define a vector field $X_t$ on $M$ by the relation
$$X_t(q)= R_{\theta(q,t)} X(q).$$
Then $X_1$ is a smooth vector field tangent to $\Zz$ at every point of 
$\Zz$. Notice that  
$\{q\in M\mid X_1(q)=0\}=\{q\in M\mid X(q)=0\}\subset \Zz$. Moreover, under 
the generic assumption that
the zeroes of $X$ are non-degenerate, the same is true for those of 
$X_1$.

Consider now the manifold  $\hat M^{+}$ obtained by gluing smoothly two 
copies of $M^+$ along $\Zz$. Since $X_1$ is tangent to $\Zz$, the vector 
field $\hat X_1$, obtained as a double copy of  $X_1|_{M^+}$, is well 
defined, continuous, and has isolated zeroes. Thus, $\chi(\hat M^+)$ is equal to the sum of the 
indices of the zeroes of $\hat X_1$. Notice that $\hat X_1(q)=0$ if and only 
if $q\in \Zz$ and $X_1(q)=0$. 
Moreover, the index of $\hat X_1$ at $q$ is equal to   
that of $X_1$, since the latter is non-degenerate. 
The same reasoning on $M^-$ and $X_1|_{M^-}$  shows that the Euler 
characteristic 
of $\hat M^-$, obtained by gluing two copies of $M^-$ along $\Zz$, is 
again equal to the sum of the indices of the zeroes of $X_1$, i.e., to $\chi(M)$. Therefore,
$$
\chi(M^+)=\frac{\chi(\hat M^+)}2=\frac{\chi(M)}2=\frac{\chi(\hat 
M^-)}2=\chi(M^-).\eqno\mbox{\quadp}
$$


\subsection{Construction of trivializable 2-ARSs with no tangency points}
\label{s-sketch}
In this section we prove Lemma~\ref{l-existence}, by showing how to construct 
a trivializable 2-ARS with no tangency 
points on every compact orientable 
two-dimensional manifold.

For the torus, an example of such structure is provided by the standard
Riemannian one. The case of a 
connected sum of two tori can be 
treated by
gluing together two copies of the pair of vector fields $X$ and $Y$ 
represented in Figure \ref{f-torichesibucano}A, which are defined on
a torus with a hole cut out. In the figure the torus is represented as a square with the standard 
identifications on the boundary. The vector fields $X$ and $Y$ are parallel on the boundary of the disk which has been cut out. 
Each vector field has exactly two zeros and the distribution spanned by $X$ and $Y$ is transversal to the singular locus.
Examples on the connected sum of three or more tori can be constructed similarly by induction. 
The resulting singular locus is represented in Figure 
\ref{f-torichesibucano}B. 
\ppotR{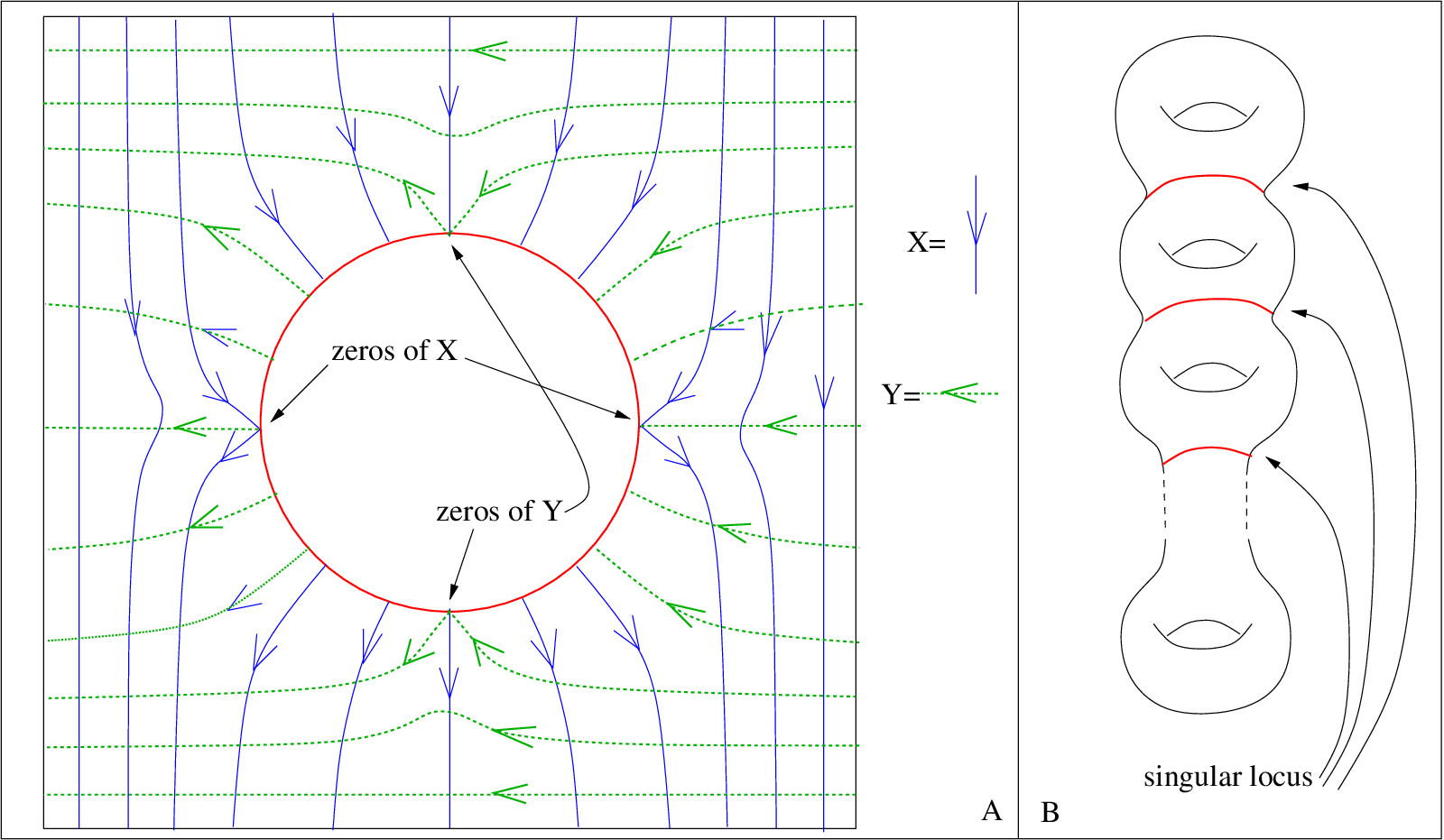}{}{17}

We are left to 
check the existence of a trivializable 2-ARS with no tangency 
points on a sphere. A simple example can be found in the literature and 
arises from a model of control of quantum systems (see \cite{q1,q4}). 
Let $M$ be a sphere in $\R^3$ centered at the origin and take 
$X(x,y,z)=(y,-x,0)$, $Y(x,y,z)=(0,z,-y)$ as orthonormal frame.
Then $X$ (respectively, $Y$) is an infinitesimal rotation around the 
third (respectively, first) axis. The singular locus is therefore given 
by the intersection of the sphere with the plane $\{y=0\}$ and none of its 
points is tangency (see Figure \ref{f-sfera}). Notice that the generic conditions given in 
Proposition \ref{p-generic} are satisfied. 

\ppotR{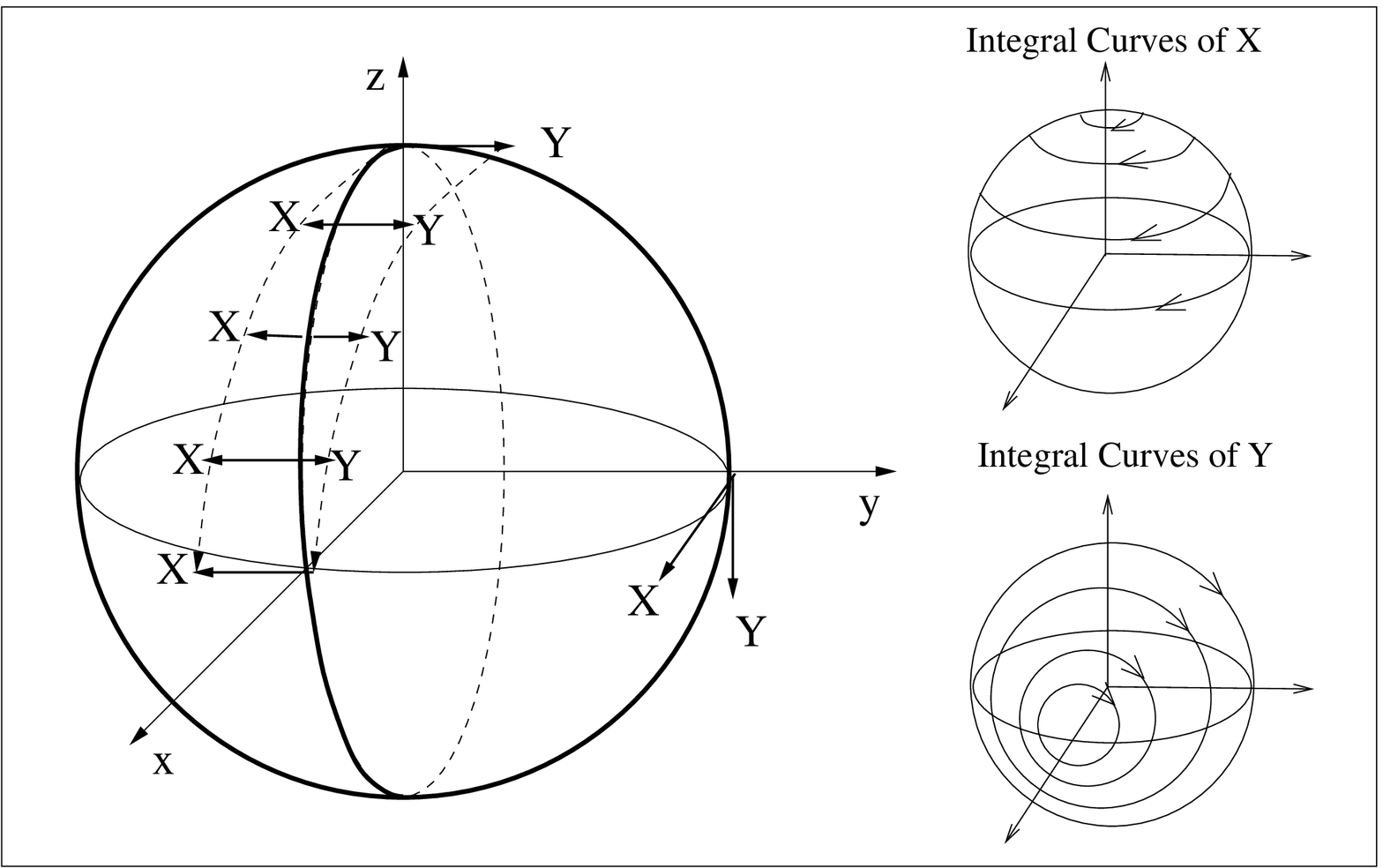}{}{12}

\end{document}